\documentclass[journal]{IEEEtran}
\pdfoutput=1
\usepackage[T1]{fontenc}
\usepackage[cmex10]{amsmath}
\usepackage{amssymb}
\usepackage{graphicx}
\usepackage{epstopdf}
\usepackage{cite}
\usepackage{url}
\usepackage{color}
\definecolor{gray}{RGB}{128,128,128}

\newcommand{\mb}{\mbox}
\newcommand{\bm}{\boldmath}
\def\b#1{\mb{\bm$#1$}}

\newlength{\fwidth}

\begin{document}

\title{Hierarchical Control Framework for Integrated Coordination between DERs and Demand Response}
\author{Di~Wu, Jianming~Lian, Yannan~Sun, Tao~Yang, and~Jacob~Hansen%
\thanks{This work was supported by the Laboratory Directed Research and Development program at Pacific Northwest National Laboratory. Pacific Northwest National Laboratory is operated for the U.S. Department of Energy by Battelle Memorial Institute under Contract DE-AC05-76RL01830.}
\thanks{D.~Wu, J.~Lian, Y.~Sun, and J.~Hansen are with the Pacific Northwest National Laboratory, Richland, WA 99354 USA (\mbox{e-mail:} \{di.wu, jianming.lian, yannan.sun, jacob.hansen\}@pnnl.gov).}
\thanks{T.~Yang is the Department of Electrical Engineering, University of North Texas, Denton, TX 76203 USA  (\mbox{e-mail:} Tao.Yang@unt.edu).}}

\maketitle

\begin{abstract}
Demand response represents a significant but largely untapped resource that can greatly enhance the flexibility and reliability of power systems.
This paper proposes a hierarchical control framework to facilitate the integrated coordination between distributed energy resources and demand response.
The proposed framework consists of coordination and device layers.
In the coordination layer, various resource aggregations are optimally coordinated in a distributed manner to achieve the system-level objectives.
In the device layer, individual resources are controlled in real time to follow the optimal power dispatch signals received from the coordination layer.
For practical applications, a method is presented to determine the utility functions of controllable loads by accounting for the real-time load dynamics and the preferences of individual customers.
The effectiveness of the proposed framework is validated by detailed simulation studies.
\end{abstract}


\section{Introduction}\label{sec:intro}
With growing emphasis on system efficiency and reliability, a great effort has been made in developing distributed energy resources~(DERs) such as distributed generator~(DG) and energy storage.
These resources are small and highly flexible compared with conventional generators, and are playing an increasingly important role in the future smart grid~\cite{Dries08,Huang08}.
On the other hand, demand-side control has presented a novel and viable way to supplement conventional supply-side control~\cite{Ipakchi09,Palen11,Callaway11}.
In fact, demand response~(DR) represents a significant but largely untapped resource in the power grid.
According to National Energy Technology Laboratory, with only 10\% customer participation, the potential nationwide value of demand dispatch could be several billion dollars per year in reduced energy costs~\cite{Goell11}.
The deployment of DERs and DR will not only defer infrastructure investments in the power grid, but also meet additional reserve requirements from renewable generation.
Although the deployment of DR and DERs can lead to more economic and reliable system operation, it requires proper coordination between DERs and DR to harvest their potential benefits.

\begin{figure*}
\centering
\includegraphics[width=0.8\textwidth]{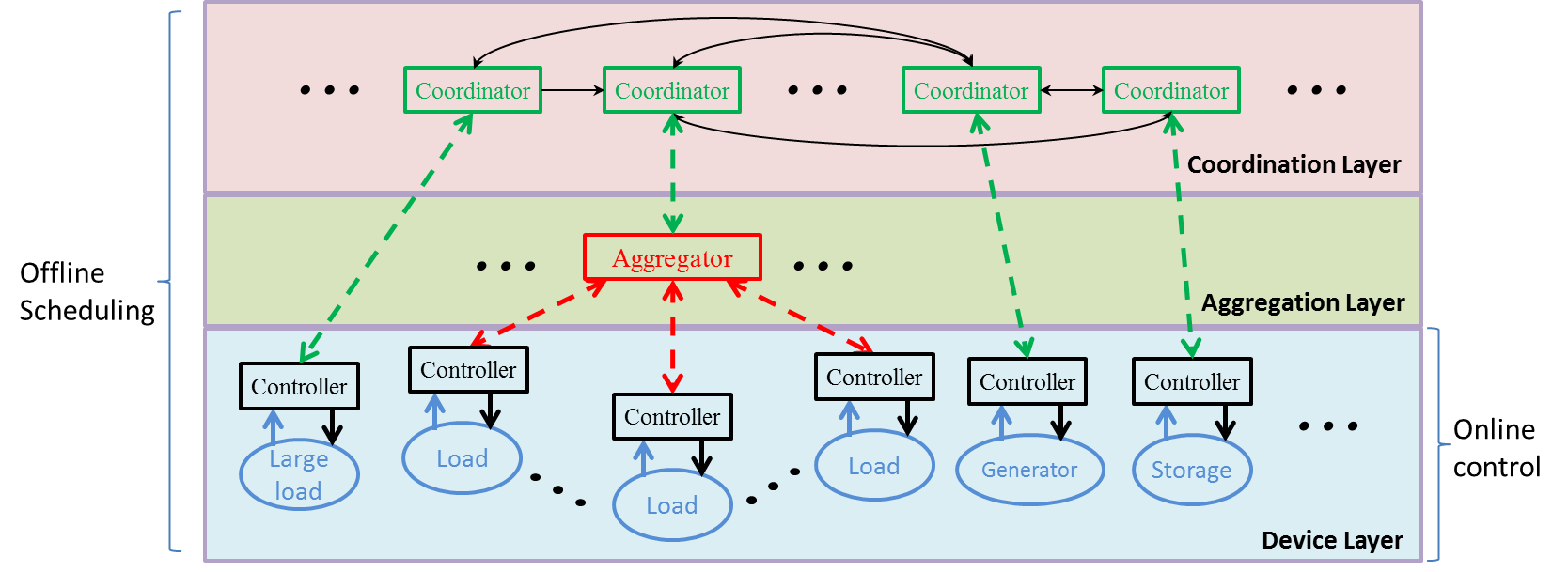}
\caption{Illustration of the proposed hierarchical control framework for integrated coordination between DERs and DR.}
\label{fig:architecture}
\end{figure*}

The coordination problem can be solved in a completely centralized manner, where a single control center accesses states of potentially thousands of devices and broadcasts control signals to them.
Such a centralized control strategy is often subject to several disadvantages, such as
high requirement and cost in communication,
substantial computational burden,
limited flexibility and scalability,
and disrespect of privacy~\cite{Farha10}.
As an alternative, a distributed control strategy has been proposed, where each control agent maintains a set of variables and updates them through information exchange with a few neighboring agents.
During the past few years, many studies have been dedicated to distributed approaches for DER coordination.
In~\cite{Domin12}, the authors developed a distributed algorithm that is resilient against potential packet drops and applied the algorithm to DER coordination.
In~\cite{Panto14}, a strategy based on the local replicator equation was presented for economic dispatch of DGs.
Other algorithms that can be applied to the DER coordination include the leader-follower consensus algorithm~\cite{Zhang11}, two-level incremental cost consensus algorithm~\cite{ZhangZ12_Convergence},
distributed algorithm based on the consensus and bisection method~\cite{XingH15}, and minimum-time consensus algorithm~\cite{YangT16}, just to name a few.
There are also studies that incorporated power losses into the distributed algorithm design~\cite{Binet14a,Elsay15}.
Recently, coordination between DERs and DR has been reported in~\cite{Hug15} and~\cite{Xu15}.
Although useful insights regarding DER and DR coordination have been reported in these studies, the existing results cannot be directly extended and applied to practical applications.
This is because the controllable loads were simply modeled as a ``generator'' with negative generation, where the load characteristics and dynamics were totally ignored. Furthermore, the studies did not address the issue of designing real-time load control strategies to achieve optimal power consumption.
This paper proposes a hierarchical control framework with two layers to achieve integrated coordination of DERs and DR.
The underlying control strategy accounts for the detailed characteristics and dynamics of controllable loads, and addresses the issue of designing real-time control strategies.

The rest of this paper is organized as follows.
In Section~\ref{sec:framework}, the major challenges of integrated coordination between DERs and DR are first discussed in detail, and then the proposed hierarchical control framework is briefly introduced with main contributions highlighted.
The top layer of the proposed framework is described in Section~\ref{sec:coordination}, where a general coordination problem between DERs and DR is formulated and solved using a distributed approach.
The bottom layer of the proposed framework is described in Section~\ref{sec:device}, where the device aggregation and real-time control are presented and illustrated using air conditioners (ACs).
In Section~\ref{sec:casestudy}, various case studies along with detailed simulation results are provided to demonstrate the effectiveness of the proposed framework.
Finally, concluding remarks are given in Section~\ref{sec:conclusion}.

\section{Problem Statement and Proposed Framework}\label{sec:framework}
Power system operation requires instantaneous power balance between generation and demand that is constantly varying.
Most balancing is achieved through energy scheduling.
In this paper, the short-term scheduling and operation problems are considered for DGs and controllable loads.
At the scheduling stage, the optimal resource allocation problem is formulated and solved between DGs and DR, where the real-time dynamics of controllable loads must be captured.
At the operation stage, real-time control is carried out so that DGs and controllable loads follow optimal power generation and consumption, respectively.

\subsection{Technical Challenges}
Although many results regarding DER and DR coordination have been presented in the literature,
there are still several technical gaps that are significant enough to prohibit practical application of these existing results.

First, the cost functions of DGs and the utility functions of controllable loads are required to formulate and solve the optimal coordination problem.
Existing studies such as~\cite{Hug15} and~\cite{Xu15} assume that those functions are available and can be directly used in the proposed distributed approaches.
However, it is not straightforward to construct the utility functions of power for controllable loads as the cost functions of power for DGs.
For instance, the utility of using an AC is directly related to the comfort an individual customer perceives at different indoor air temperatures rather than the power consumption.
Therefore, it is required in practice to extract the utility functions and capture the underlying economics based on the preferences of individual customers.

Second, it is required for practical applications to consider the operation stage as well.
After the coordination problem is solved at the scheduling stage, individual resources are expected to follow optimal generation or consumption through real-time control.
It is straightforward for DGs to meet this expectation because their generation level can be continuously adjusted with existing generator controllers.
However, this is often not the case for controllable loads.
Some controllable loads such as thermostatically controllable loads have not been designed with the capability to continuously adjust their power consumption.
Furthermore, their power consumption cannot be directly controlled and is usually indirectly affected by other control variables.
For example, the thermostat of an AC receives the temperature setpoint as the control input and then automatically switches the compressor on and off to maintain the indoor air temperature around the setpoint.
Therefore, a real-time load controller has to be designed for individual controllable loads using the locally acceptable control input while capturing the underlying economics.

Effectively coordinating and controlling DERs and DR for short-term scheduling and real-time control cannot be realized by simply adding one coordination algorithm to another load control approach.
A systematic method is needed to capture the underlying economics and dynamics of controllable load synthetically in both scheduling and real-time control.
The proposed framework herein exactly meets such a need.
The main contributions of this paper are summarized as follows:
\begin{itemize}
      \item We discuss the gap between short-term scheduling and real-time control from existing coordination algorithms and load control approaches, and identify challenges to bridging the gap.
      \item We design a holistic hierarchical control framework, which is capable of directly adopting existing coordination algorithms and load modeling/control approaches.
      \item We define the functionality and formulate mathematic problems in each layer, and specify the information to be exchanged between layers/sublayers in both short-term scheduling and real-time control.
      \item We identify candidate coordination algorithms and load control approaches that can fit in the proposed framework, and use example algorithm and approach to illustrate the proposed framework.
\end{itemize}

\subsection{Proposed Framework}
To overcome the technical challenges described above, this paper proposes a hierarchical control framework as shown in Fig.~\ref{fig:architecture} to facilitate integrated coordination between DERs and DR.
The proposed framework consists of two layers including coordination (top) and device layers (bottom).
In Fig.~\ref{fig:architecture}, dash lines represent information flow between layers/sublayers, where the information exchange frequency is the same as short-term scheduling.
Solid lines represent information flow within a layer/sublayer, where the information exchange frequency is typically much higher than short-term scheduling.
An overview of each layer is provided herein and more details are provided in the following sections.

\begin{itemize}
 \item The coordination layer is only involved in short-term scheduling stage.
 Prior to each scheduling period, each coordinator receives aggregated utility or cost functions from aggregators or device controllers, as indicated by the green dash lines with up arrows.
 Then, the aggregation of various resources including DGs and controllable loads are optimally coordinated to achieve power balance.
 To overcome the disadvantages associated with centralized coordination algorithms, a distributed coordination method can be employed, where local variables are exchanged iteratively following algorithms as explained in Section III.
 Once the coordination problem is solved, the regulation signals are sent back to aggregators or device controllers for real-time control, as indicated by the green dash lines with down arrows.
  \item The device layer includes two sublayers: device aggregation and device control.
 \begin{itemize}
\item In the aggregation sublayer, DERs are divided into groups as appropriate.
Prior to each scheduling period, each aggregator received utility or cost functions from its underlying device controllers (as indicated by the red dash lines with up arrows), determines the aggregated functions,
and then sends these information to the corresponding coordinator in the top layer (as indicated by the green dash lines with up arrows).
After coordination problem is solved, each aggregator receives the regulation signals from top layer (as indicated by the green dash lines with down arrows) and then broadcasts these signals to its underlying devices (as indicated by the red dash lines with down arrows) to collectively provide the desired power generation or consumption.
\item The device control sublayer is involved in both scheduling and real-time operation stages.
At the scheduling stage, the controller at each device reports to its commander (either an aggregator or a coordinator) the required information (as indicated by the green/red dash lines with up arrows), which is then used for coordination.
After coordination problem is solved, it receives the regulation signals from its commander (as indicated by the green/red dash lines with down arrows).
During real-time operation, it regulates devices to fulfill their functionality (e.g., control the indoor air temperature within the comfort zone) while following the scheduled energy generation or consumption.
\end{itemize}
\end{itemize}

\section{Optimal Resource Coordination}\label{sec:coordination}

\subsection{General Description}
In the coordination layer, the optimal coordination problem is solved for each coordination period to determine the most economic schedule of power generation and consumption for DERs and DR, respectively.
The objective is to maximize the social welfare, i.e., the difference between the utility of power consumption and the cost of power generation, while meeting the desired total power output without violating operating constraints of individual resources.

The mathematical formulation of the scheduling problem for each coordination period is presented as follows,
where power generation or consumption should be understood as the average value during each coordination period,
\begin{IEEEeqnarray}{r'c}
\label{eq:primal}
\min_{p_{G_i},\, p_{L_j}} & \sum_{i=1}^{N_G} C_{G_i}(p_{G_i}) - \sum_{j=1}^{N_L} U_{L_j}(p_{L_j})\IEEEyesnumber \IEEEyessubnumber \label{eq:primala}\\
\text{subject to} & \sum_{i=1}^{N_G} p_{G_i}-\sum_{j=1}^{N_L} p_{L_j}=D \IEEEyessubnumber \label{eq:primalb}\\
& 0\le P^{\min}_{G_i} \le p_{G_i} \le P^{\max}_{G_i} \IEEEyessubnumber \label{eq:primalc}\\
& 0\le P^{\min}_{L_j} \le p_{L_j} \le P^{\max}_{L_j} \IEEEyessubnumber \label{eq:primald}
\end{IEEEeqnarray}
where different notations are defined as follows:
\begin{itemize}
\renewcommand\labelitemi{--}
  \item $N_G$ (or $N_L$) is the number of generator (or load) aggregation in the network;
  \item $C_{G_i}(p_{G_i})$ is the cost of the $i$-th generator aggregation as a function of the power generation $p_{G_i}$;
  \item $U_{L_j}(p_{L_j})$ is the utility of the $j$-th load aggregation as a function of the power consumption $p_{L_j}$;
  \item $D$ is the desired total power output;
  \item $[P^{\min}_{G_i},P^{\max}_{G_i}]$ is the range of power generation for the $i$-th generator aggregation;
  \item $[P_{L_j}^{\min},P^{\max}_{L_j}]$ is the range of power consumption for the $j$-th load aggregation.
\end{itemize}
Note that $[P_{L_j}^{\min},P^{\max}_{L_j}]$ for each load aggregation can be obtained by aggregating the power range of individual controllable loads.
For example, the average power consumption of an AC for the next 5~minutes depends on the temperature setpoint selected by the homeowner,
the current indoor air temperature and the outside air temperature, etc.

\subsection{Proposed Approach}
Without loss of generality, all the DGs can be enumerated as the first $N_G$ agents.
The optimization problem defined in~\eqref{eq:primala}--\eqref{eq:primald} can then be generalized as
\begin{IEEEeqnarray}{r'c'l}
\label{eq:primal1}
\min_{p_{i}} &  \sum_{i=1}^N C_i(p_i) &\IEEEyesnumber \IEEEyessubnumber \label{eq:primal1_cost}\\
\text{subject to} & \sum_{i=1}^N p_i=D& \IEEEyessubnumber\\
& P_i^{\min} \le p_i \le P_i^{\max}, & i=1,\,\ldots,\,N \IEEEyessubnumber
\end{IEEEeqnarray}
where $N=N_G+N_L$,
\begin{equation}
p_i = \left\{ \begin{array}{ll}
p_{G_i}, & i=1,\ldots,N_G\\
-p_{L_{i-N_G}}, & i=N_G+1,\ldots,N_G+N_L
\end{array}\right.
\end{equation}
\begin{equation}
C_i(\cdot) = \left\{ \,
\begin{IEEEeqnarraybox}[][c]{l?l}
\IEEEstrut
C_{G_i}(\cdot), &i=1,\ldots,N_G\\
-U_{L_{i-N_G}}(\cdot), &i=N_G+1,\ldots,N_G+N_L \label{eq:utilitycostfunction}.
\IEEEstrut
\end{IEEEeqnarraybox}
\right.
\end{equation}
The optimal solution can be obtained through various distributed coordination algorithms reviewed in Section~\ref{sec:intro}.
Most of these algorithms are consensus-based with marginal cost modeled as consensus variables.
They solve the problem essentially through price-directive decomposition, which is actually the gradient method applied to the dual problem. 
Different methods have been proposed to update the dual variable using partial or total mismatch between demand and supply (the gradient of the dual problem)\cite{Zhang11,ZhangZ12_Convergence,Domin12,Hug15,XingH15}.
In this paper, we use the distributed coordination algorithm proposed in~\cite{Kar12}.

Prior to each scheduling period, each coordinator receives the aggregated utility or cost functions as well as the power operating ranges from aggregators or device controllers, as indicated by the green dash lines with up arrows in Fig.~\ref{fig:architecture}.
Next, each coordinator converts the received cost/utility functions and power output to $C_i(\cdot)$ and $p_i$, respectively, according to \eqref{eq:primal1}.
Then, the coordinator starts to run the coordination algorithm as shown in \eqref{distributed-update-LL}.
\begin{IEEEeqnarray}{rcl}\label{distributed-update-LL}
\lambda_{i}(k+1)&=& \lambda_i(k)-\beta_k \sum_{j\in \mathcal{N}_i} (\lambda_i(k)-\lambda_j(k)) \nonumber \\
&&\hspace{50pt} -\alpha_k (p_{i}(k)-D_i) \,\text{,}  \IEEEyesnumber \IEEEyessubnumber  \label{distributed-update1-LL}\\
p_{i}(k+1)&=&\nabla C^{-1}_i(\lambda_{i}(k+1))\,\text{,}  \IEEEyessubnumber \label{distributed-update2-LL}
\end{IEEEeqnarray}
where $\mathcal{N}_i=\{j \in \mathcal{V}|(j,i)\in \mathcal{E}\}$ is the neighboring set of the $i$-th agent, $\nabla C_i(\cdot)$ is the derivative of cost function and $\nabla C^{-1}$ denotes its inverse function,
$\alpha_k$ and $\beta_k$ are the gain parameters at step $k$ for innovation term and consensus term , respectively,
and $D_i$ is chosen such that $\sum_{i=1}^{N}D_i=D$.
The determination of $D_i$ can be arbitrary.
In practice, one option to determine $D_i$ is that the system operator forecasts the total demand $D$, and then arbitrarily distributes this demand to a small set of agents or even a single agent ($D_i$ is zero for the remaining agents).
Alternatively, each agent can also determine its own $D_i$.
This strategy is used in this paper. Please refer to Section V.A for more details.

With this algorithm, each coordinator $i$ only maintains a local variable $\lambda_{i}$ that is the estimate of the optimal incremental cost,
and updates it through information exchange with its neighboring coordinators (as indicated by the black lines in the coordinator layer in Fig.~\ref{fig:architecture}.)
By executing~\eqref{distributed-update-LL}, $\lambda_i(k)$ and $p_i(k)$ at each coordinating agent will converge to the optimal dual variable (clearing prices) and power output,
which are sent back to aggregators or device controllers for real-time control,
as indicated by the green dash lines with down arrows in Fig.~\ref{fig:architecture}.

\section{Device Aggregation and Control}\label{sec:device}

\subsection{General Description}
Although it is necessary to have $C_{G_i}(p_{G_i})$ and $U_{L_j}(p_{L_j})$ to formulate the optimal coordination problem as shown in~\eqref{eq:primala}--\eqref{eq:primald}, the distributed algorithms only require their derivatives, $C^\prime_{G_i}(p_{G_i})$ and $U^\prime_{L_j}(p_{L_j})$, to solve this coordination problem.
The derivative of $C_{G_i}(p_{G_i})$ (or $U_{L_j}(p_{L_j})$) is often referred to as the supply (or demand) curve, which characterizes the relationship between marginal cost (or utility) and power generation (or consumption).
Hence, the device layer has to determine the supply and demand curves of various resources and send them to the coordination layer at the beginning of each coordination period.
To maintain the scalability of the proposed framework in dealing with a large number of resources, the device layer is further divided into two sublayers: device aggregation and device control.
In the device aggregation sublayer, resources of similar type are grouped together when their individual sizes are small.
In practice, the resource aggregation is usually employed to either facilitate coordination processes or represent business models.
Each aggregator serves as the message channel between the coordination layer and the device control sublayer.
It collects the individual supply or demand curves from resources within its aggregation (as indicated by the red dash lines with up arrows in Fig.~\ref{fig:architecture}),
and then sends the aggregated curve to the coordination layer (as indicated by the green dash lines with up arrows in Fig.~\ref{fig:architecture}).
Then it sends the optimal dispatch signals received from the coordination layer to the device control sublayer (as indicated by the red dash lines with down arrows in Fig.~\ref{fig:architecture}).
In the device control sublayer, real-time control translates the optimal dispatch signals into local control inputs so that individual resources can follow the optimal generation or consumption.

The supply curves of DGs can be easily determined based on generator operational cost, fuel efficiency, and fuel cost.
It is also straightforward for them to follow the optimal generation in real time because their generation level can be continuously adjusted with well established controllers.
However, for controllable load, technical challenges exist in
i) determining the demand curves of controllable loads based on individual customer preference, and ii) designing real-time control that can translate optimal power demand into locally acceptable control input.
As pointed out in Section~\ref{sec:framework}, one essential step is to obtain the relationship between marginal utility and local control input.

\subsection{Proposed Approach}\label{sec:demand-response-curve}
The demand curve dynamically represents how individual customers value convenience or comfort and the corresponding energy usage.
It is essential to capture the opportunity cost of DR.
To extract the demand curve, it is necessary to quantitatively relate the marginal utility of individual customers for power demand to the local control input based on customer preference.
Herein, a practical method is presented for ACs to extract such a relationship.
This method was originally proposed in the GridWise\textsuperscript{\textregistered} demonstration project~\cite{Fuller11}, and then rigorously analyzed in~\cite{LiS15a}.
Although it has been specifically presented for ACs, the underlying control philosophy can be easily extended and applied to other types of controllable loads.

\begin{figure}
\centering
\includegraphics[width=0.45\textwidth]{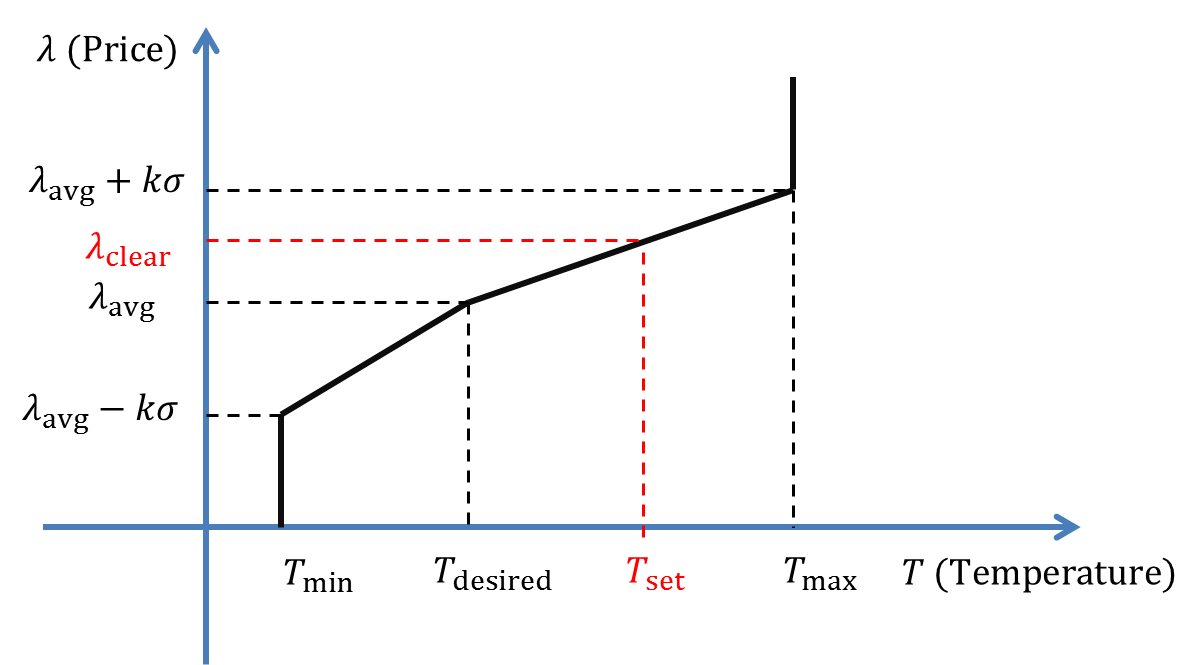}
\caption{Illustration of the response curves for air conditioners.}
\label{fig:response}
\end{figure}

\begin{figure}
\centering
\includegraphics[width=0.35\textwidth]{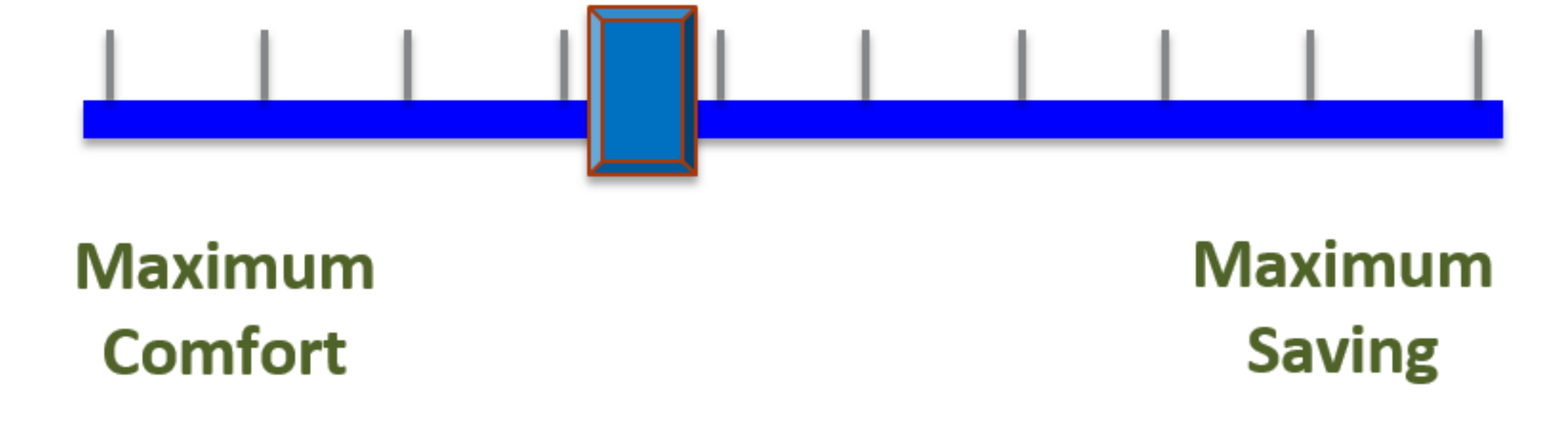}
\caption{User interface in the GridWise\textsuperscript{\textregistered} demonstration project}
\label{fig:interface}
\end{figure}

This method represents the relationship between marginal utility and local control input by a response curve as illustrated in Fig.~\ref{fig:response}, which is determined by several parameters.
The parameters $\lambda_\text{avg}$ and $\sigma$ are the average and variance, respectively, of the electricity prices over a period of time in the past, which can be calculated by a local controller or load aggregator.
The parameters $T_\text{desired}$, $T_\text{min}$, and $T_\text{max}$ are directly specified by users, where $T_\text{desired}$ is the desired indoor air temperature setpoint, and $T_\text{min}$ and $T_\text{max}$ are the lower and upper bounds of the acceptable indoor air temperature setpoint.
The parameter $k$ is a positive number completely abstracted from the owner's preference of indoor air temperature setpoint over the electricity price. For example, when $k$ is very large, the response curve becomes an almost vertical line at $T_\text{desired}$.
This implies that the homeowner is very sensitive to the indoor air temperature, and would like to maintain the indoor air temperature setpoint at $T_\text{desired}$ regardless of the electricity price.
When $k$ is close to zero, the response curve becomes an almost horizontal line at $\lambda_\text{avg}$.
This implies that the house owner is very sensitive to the electricity price, and would like to sacrifice comfort for cost saving.
In the GridWise\textsuperscript{\textregistered} demonstration project, the abstraction of $k$ is done by letting individual homeowners specify their preferences of comfort over cost through a user interface as shown in Fig.~\ref{fig:interface}.

\begin{figure}
\centering
\includegraphics[width=0.45\textwidth]{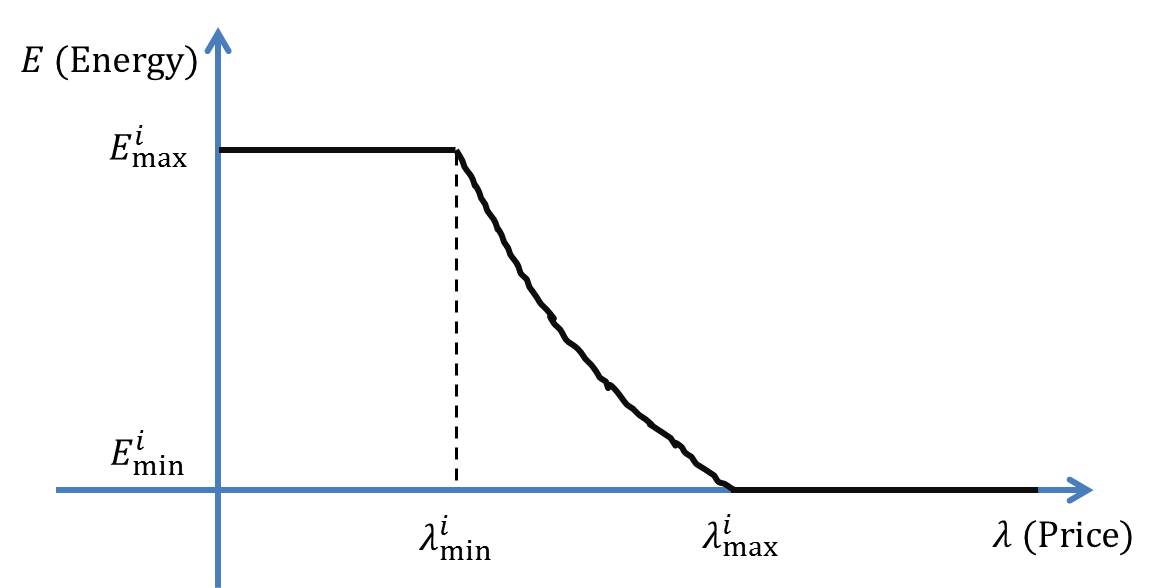}
\caption{Illustration of the demand response curves of air conditioners.}
\label{fig:ACdemand}
\end{figure}

In the proposed framework, prior to each coordination period, each load controller at the device layer determines its DR curve (which is equivalent to the utility/cost function) considering the load dynamics.
These curves are sent to aggregators (as indicated by the red dash lines with up arrows in Fig.~\ref{fig:architecture}), and
then aggregated at aggregators and later used by coordinators to run the distributed coordination algorithms for determining the optimal power consumption and real-time price signal.
The price signal is then broadcasted to each device every 5 minutes, as indicated by the red dash lines with down arrows in Fig.~\ref{fig:architecture}.
Within each 5-minute operating period, this response curve will be used by a real-time controller to translate the optimal power demand represented by clearing price $\lambda_\text{clear}$ into the indoor air temperature setpoint $T_\text{set}$ for this operating period.
Recall that the demand curve is the mapping from marginal utility to power demand.
Since the response curve is the mapping from marginal utility to the indoor air temperature setpoint, it is only left to determine the relationship between indoor air temperature setpoint and power demand.
The dynamics of each AC~$i$ can be described by the Equivalent Thermal Parameter model.
The detailed model can be found in \cite{UR,Thomas12}, and can be represented in a simplified form as
\begin{equation}
\label{eq:etp}
\dot{\b x}_i(t)=\left\{\begin{array}{ll}
\b A_i\b x_i(t)+\b B_\textrm{on}^i & \textrm{if $q_i(t)=1$} \\
\b A_i\b x_i(t)+\b B_\textrm{off}^i & \textrm{if $q_i(t)=0$}, \\
\end{array}\right.
\end{equation}
where $\b x_i(t)$ is the continuous state vector consisting of indoor air temperature $T_a^i(t)$ and mass temperature $T_m^i(t)$,
and $q_i(t)$ denotes the operating mode of the AC with $q_i(t)=1$ when it is ON and $q_i(t)=0$ when it is OFF.
The operating mode of the AC for cooling is usually controlled by a hysteretic controller,
\begin{equation}
\label{eq:hyster}
q_i(t^+)=\left\{\begin{array}{ll}
1 & \textrm{if $T_a^i(t)\ge T_\textrm{set}+\delta/2$} \\
0 & \textrm{if $T_a^i(t)\le T_\textrm{set}-\delta/2$} \\
q_i(t) & \textrm{otherwise}, \\
\end{array}\right.
\end{equation}
where $\delta$ is the hysteresis band centered around the indoor air temperature setpoint $T_\text{set}$.
When the model parameters $(\b A_i,\b B_\textrm{on}^i,\b B_\textrm{off}^i)$ are known to local controllers, the relationship between indoor air temperature setpoint and power consumption can be derived,
which finally leads to the demand curve as shown in Fig.~\ref{fig:ACdemand} by taking into account the corresponding response curve.
The determination of $E_\text{max}^i$, $E_\text{min}^i$, $\lambda_\text{min}^i$, and $\lambda_\text{max}^i$ is provided in the Appendix.
If the model parameters are unknown, they can be estimated based on the measured indoor air temperature as proposed in~\cite{LiS15b}.
The individual demand curves will be sent to the load aggregator, where they are aggregated together, and sent to the coordination layer to solve the optimal coordination problem.

With this method, the load utility/cost functions depend on the market clearing price during previous periods, the outside air temperature, and customer preference for comfort.
Therefore, these functions in \eqref{eq:utilitycostfunction} are time-varying from one scheduling period to another.
The DR curve illustrated in Fig.~\ref{fig:ACdemand} is the derivative of utility function \eqref{eq:primala}.
Since the DR curve is monotonically non-increasing, the utility function is concave.
The negative utility function in \eqref{eq:primala} becomes a cost function in \eqref{eq:primal1_cost}, which is convex.

\section{Case Studies}\label{sec:casestudy}
This section demonstrates the proposed hierarchical control framework by case studies on the IEEE 123-node system that was prepared by IEEE PES Distribution System Analysis Subcommittee's Distribution Test Feeder Working Group~\cite{testfeeders}.
The simulation studies are implemented in GridLAB-D \cite{gridlabd}, which is an advanced open-source power systems modeling and simulation environment developed at Pacific Northwest National Laboratory.

\subsection{Test System Description}
\begin{figure}
	\centering
	\includegraphics[width=0.485\textwidth]{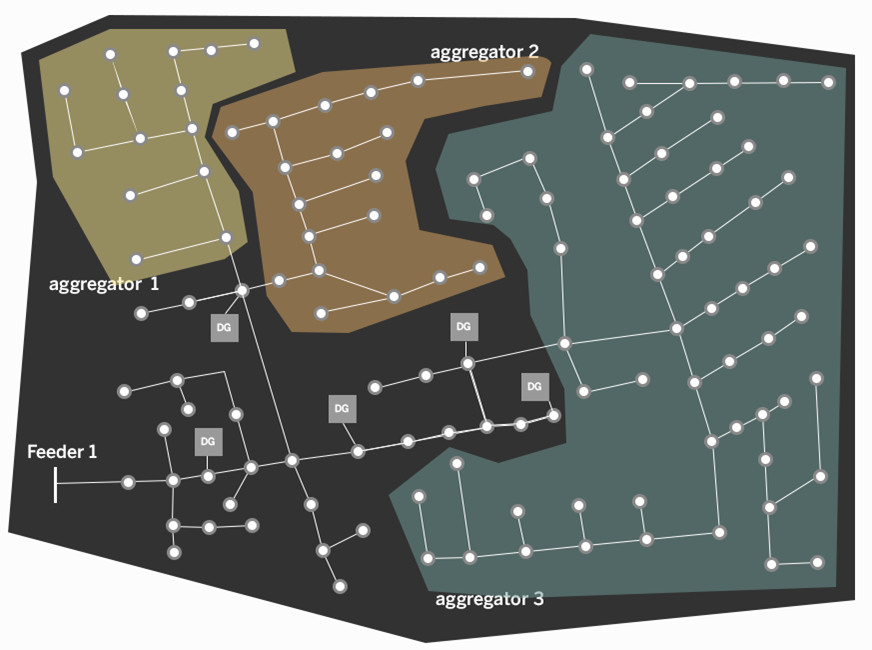}
	\caption{IEEE 123-node test system.}
	\label{fig:testCircuit}
\end{figure}

The IEEE 123-node test system shown in Fig.~\ref{fig:testCircuit} consists of 123 nodes and 118 lines.
It has been modified to include houses with ACs and other residential loads.
The number of houses has been adjusted to match the peak load provided in the test system dataset, which results in 1,222 houses.
The following studies assume that 988 ACs participate in the DR program, and the remaining 234 ACs are uncontrollable as other residential loads.
The controllable ACs are grouped into three load aggregations, where the number of houses under each aggregation are 98, 254, and 632, respectively.
There are five DGs connected to the system, whose generation cost parameters are adopted from~\cite{Logen12} and~\cite{Modir13} and listed in Table~\ref{table-DG}.
Although it is typical to represent generation cost by quadratic functions, the distributed algorithms reviewed in the Introduction section are able to handle convex functions that are more general than quadratic functions.
In fact, the cost functions of DR in this paper are convex, as explained in Section~\ref{sec:demand-response-curve}.
The distributed algorithm for solving the optimal coordination problem is selected to be the leaderless algorithm defined in~\eqref{distributed-update1-LL} and~\eqref{distributed-update2-LL}, which requires undirect communication networks.

\begin{table}[!t]
\caption{Generator Parameters}
\renewcommand{\arraystretch}{1.2}
\label{table-DG}
\centering
\begin{tabular}{|c|c|c|c|c|}
\hline
DG No. & $a_i \, \text{(kW}^2 \text{h)}$ & $b_i$ (\$/kWh)  & $c_i$ (\$/h) &  Range (kW) \\
\hline
1 & 0.00015 & 0.0267 & 0.38 & [50,500] \\
\hline
2 & 0.00052 & 0.0152 & 0.65 & [20,100] \\
\hline
3 & 0.00042 & 0.0185 & 0.4 & [40,200] \\
\hline
4 & 0.00031 & 0.0297 & 0.3 & [20,250] \\
\hline
5 & 0.00025 & 0.0156 & 0.33 & [30,300] \\
\hline
\end{tabular}
\end{table}

Prior to each coordination period, which is selected to be 5~minutes, controllable loads are aggregated and then coordinated with DGs so that they can meet
\begin{equation}
\label{eq:balance}
\text{contr. load} + \text{uncontr. load} - \text{DG gen.}=\text{ref.}
\end{equation}
In \eqref{eq:balance}, ref. represents the desired power consumption of the distribution system and can be commanded by system operators.
For example, in island microgrid operation, the command will be set to zero, i.e., DGs and controllable loads are scheduled to balance the uncontrollable load within the microgrid for each 5-minute period.
For a grid-connected distribution system, DGs and DR can be controlled to follow a given signal for either reducing energy cost or providing grid services.
For example, they can actively participate in the load following service by setting the reference signal as
$$\text{ref.}=\text{feeder hourly schedule}-\text{load following signal}.$$
In the coordination problem \eqref{eq:primal1}, the desired total output from DGs and DR is
\begin{equation}
  D=\text{ref.} -  \text{uncontr. load} \,.
\end{equation}
In order to apply the distributed algorithm in \eqref{distributed-update-LL}, we need to determine $D_i$ such that $\sum_i D_i=D$.
In this paper, we choose to distribute $D$ only among three load aggregators and set $D_i$ to zero for DGs.
The distribution of $D$ is realized by distribution of ref. signal and local uncontr. load.
First, ref. signal can be arbitrarily distributed among load aggregators offline, prior to scheduling period.
In this paper, we evenly distribute this signal among three aggregators.
On the other hand, the uncontrollable load at each aggregator is unknown and needs to be forecasted.
While there are different load forecasting methodologies, e.g., artificial neural networks~\cite{Hippe01}, autoregressive moving average models~\cite{Huang03}, and semi-parametric additive models~\cite{Fan12},
this work simply assumes that the forecast of local uncontrollable load in the next 5-minute period is equal to the measured local uncontrollable load in the current period.
Please note that renewable generation as one kind of important distributed generation could also be modeled as negative uncontrollable load as it is typically controlled for maximum power tracking.
It is forecasted and adjusted every 5 minutes, and becomes a component of `uncontr. load'.
In such a way, controllable load and DGs can be optimally used to help address the uncertainty and variability from renewable generation.

When performing DGs and DR coordination for the time period $t$ using~\eqref{distributed-update1-LL} and~\eqref{distributed-update2-LL}, the initial values of power $p_i$ and marginal costs $\lambda_i$ are set to be converged values in the time period $t-1$.
This can help to greatly reduce the required number of iterations.

\subsection{Simulation Results}
\subsubsection{Base case (Case 1)}
\begin{figure}
\centering
\includegraphics[width=0.45\textwidth]{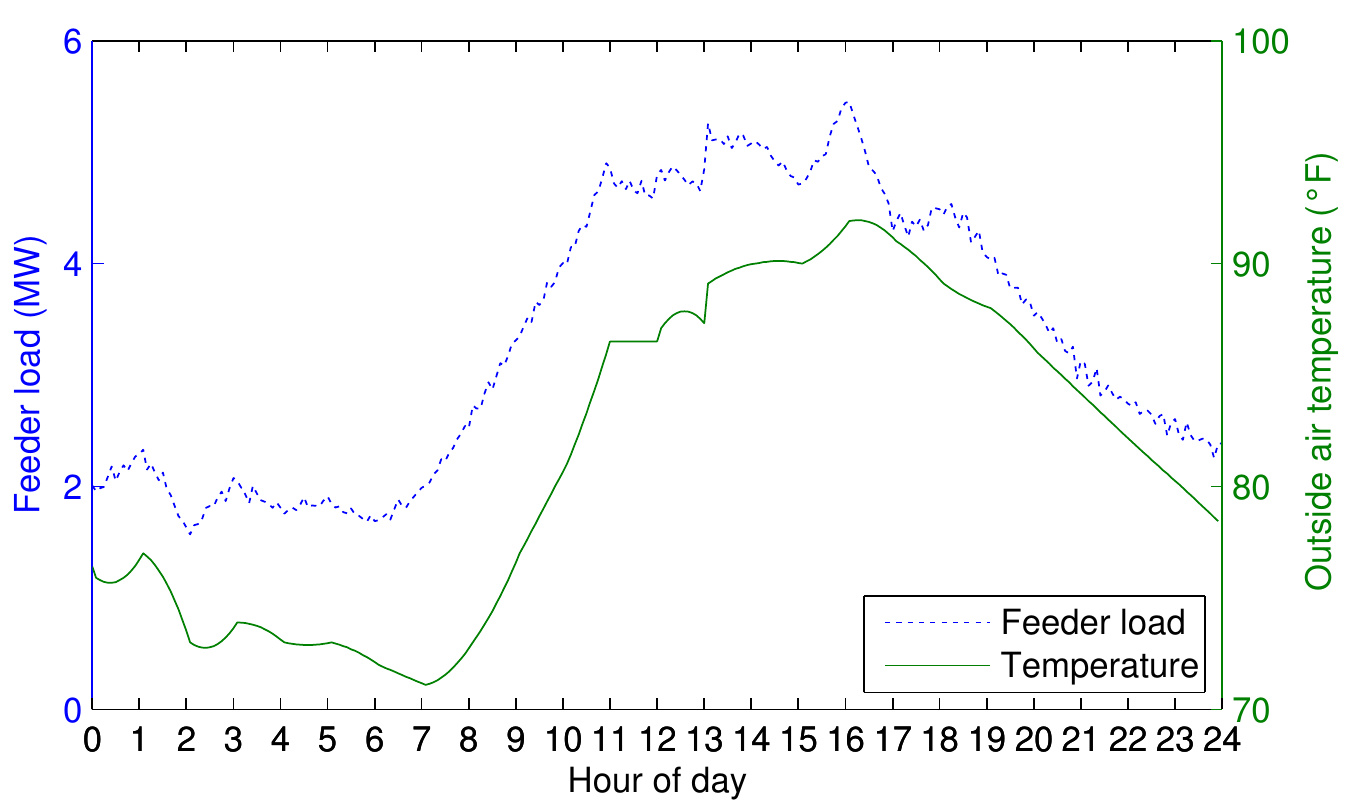}
\caption{Base case feeder load (5-minute average) and outside air temperature.}
\label{fig:base}
\end{figure}

The test system is first simulated without any DGs and controllable loads for a typical summer day with a minimum time step of 30 seconds.
The 5-minute average feeder power consumption is plotted in Fig.~\ref{fig:base}, together with the outside air temperature.
Since AC load accounts for more than 80\% of the total load in this system, the system load increases as the outside air temperature rises.

\subsubsection{Active DER and DR Coordination (Case 2)}
\begin{figure}[!t]
\centering
\includegraphics[width=0.45\textwidth]{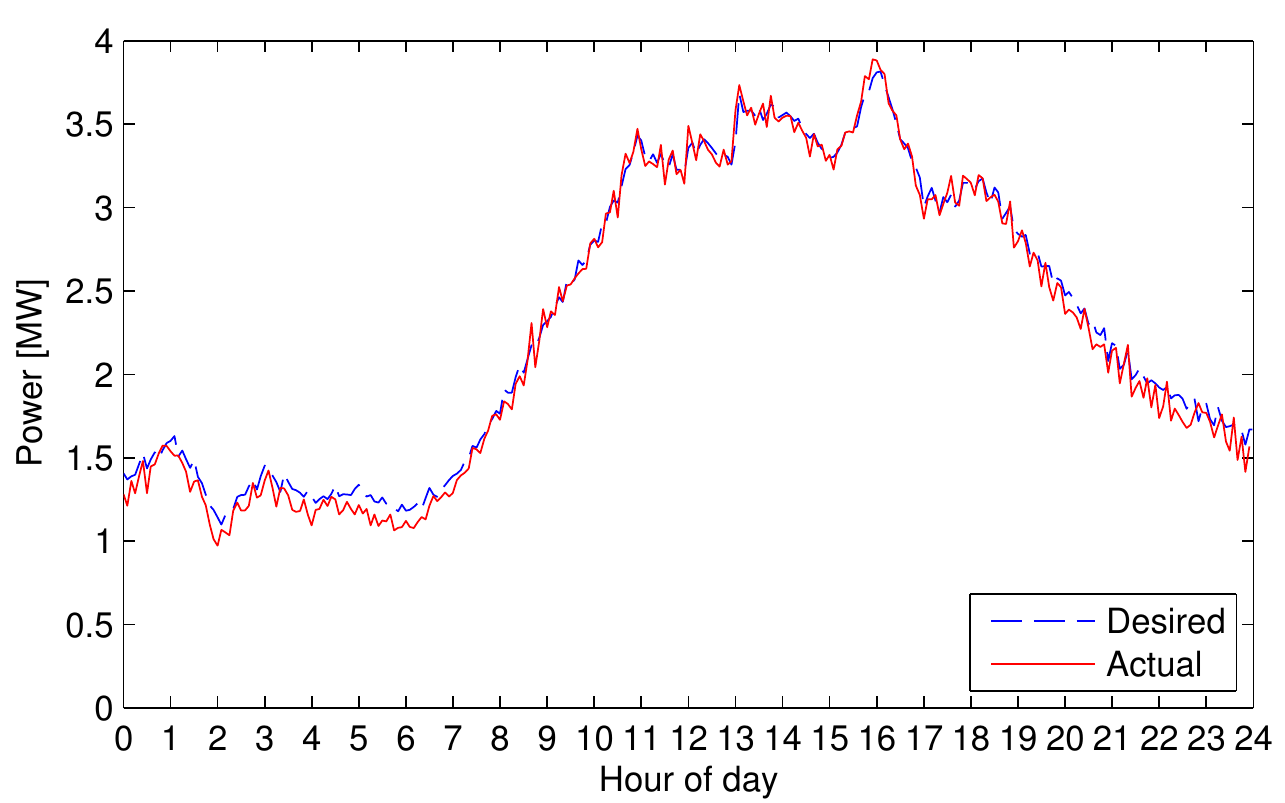}
\caption{Desired vs. resulted feeder load (both are  5-minute average).}
\label{fig:Refer}
\end{figure}

The test system is then simulated with DGs and controllable loads under the proposed hierarchical control framework for the same summer day.
In general, the reference signal can be any time series within the capability of the active distribution system.
To verify the effectiveness of the proposed framework, the desired feeder load consumption is set to be 0.7 of the feeder load in base case, as shown by the blue dashed line in Fig.~\ref{fig:Refer}.
Such a reference signal is simple to construct yet useful for testing the proposed method because
\begin{itemize}
\item The 30\% reduction of load at the feeder requires the participation of DG and DR during scheduling, which is exactly what we need to study.
\item The reduction is proportional to the load feeder in the base case and therefore varies with time.
Such varying load reduction requires DG and DR to vary their generation or consumption in a coordinative manner.
\item Such a desired signal requires DG and DR to support the local system more during peak hours than off-peak hours, which seems plausible.
The test case enables us to compare DER participation in peak hours with off-peak hours, as well as the difference in energy price of the distribution system.
\end{itemize}

The obtained 5-minute average power consumption is plotted by the red curve in Fig.~\ref{fig:Refer}.
As can be seen, the actual feeder load follows the desired value with reasonable accuracy.
The small mismatch is due to a few factors such as approximation of demand curve, errors in uncontrollable load forecast, and approximation of optimal solution in coordination layers.

\begin{figure}
\centering
\includegraphics[width=0.45\textwidth]{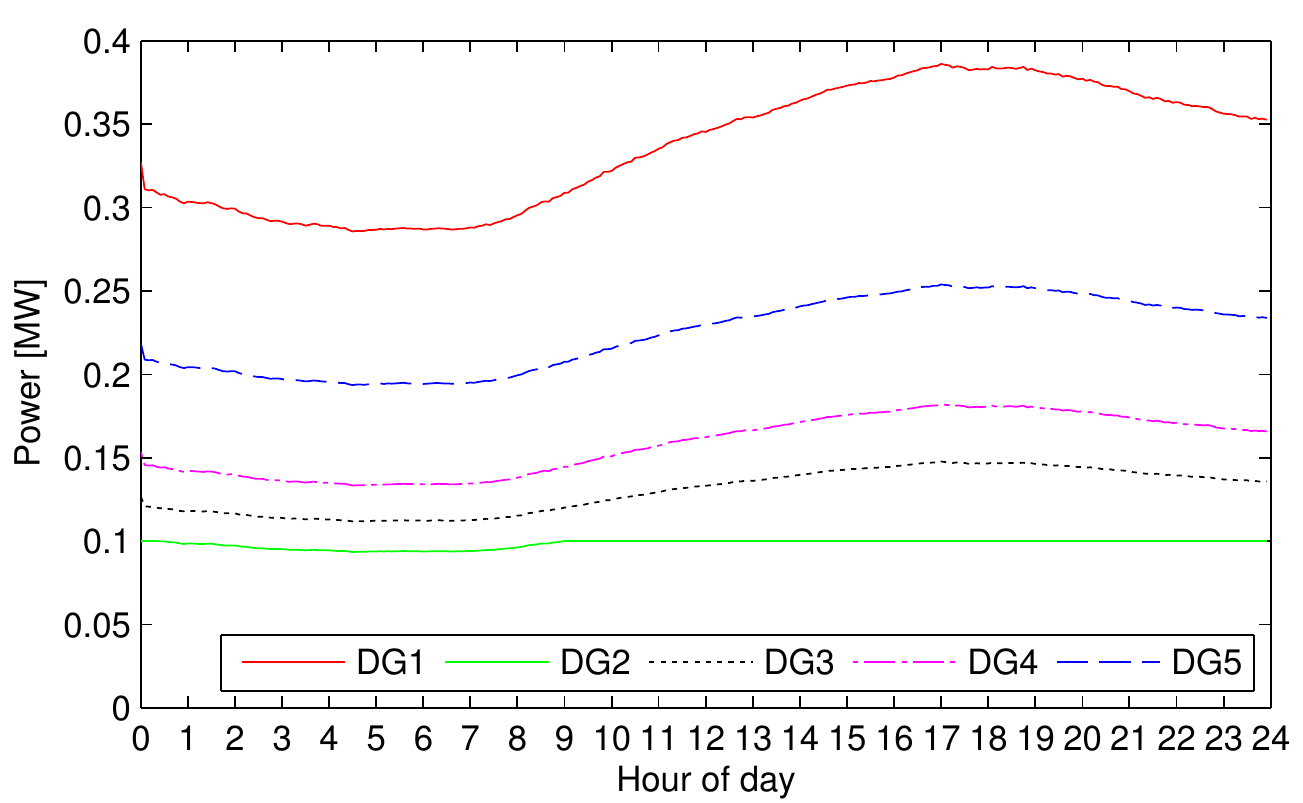}
\caption{Generation output from DGs in Case 2.}
\label{fig:DG_gen}
\end{figure}

\begin{figure}[!t]
\centering
\includegraphics[width=0.45\textwidth]{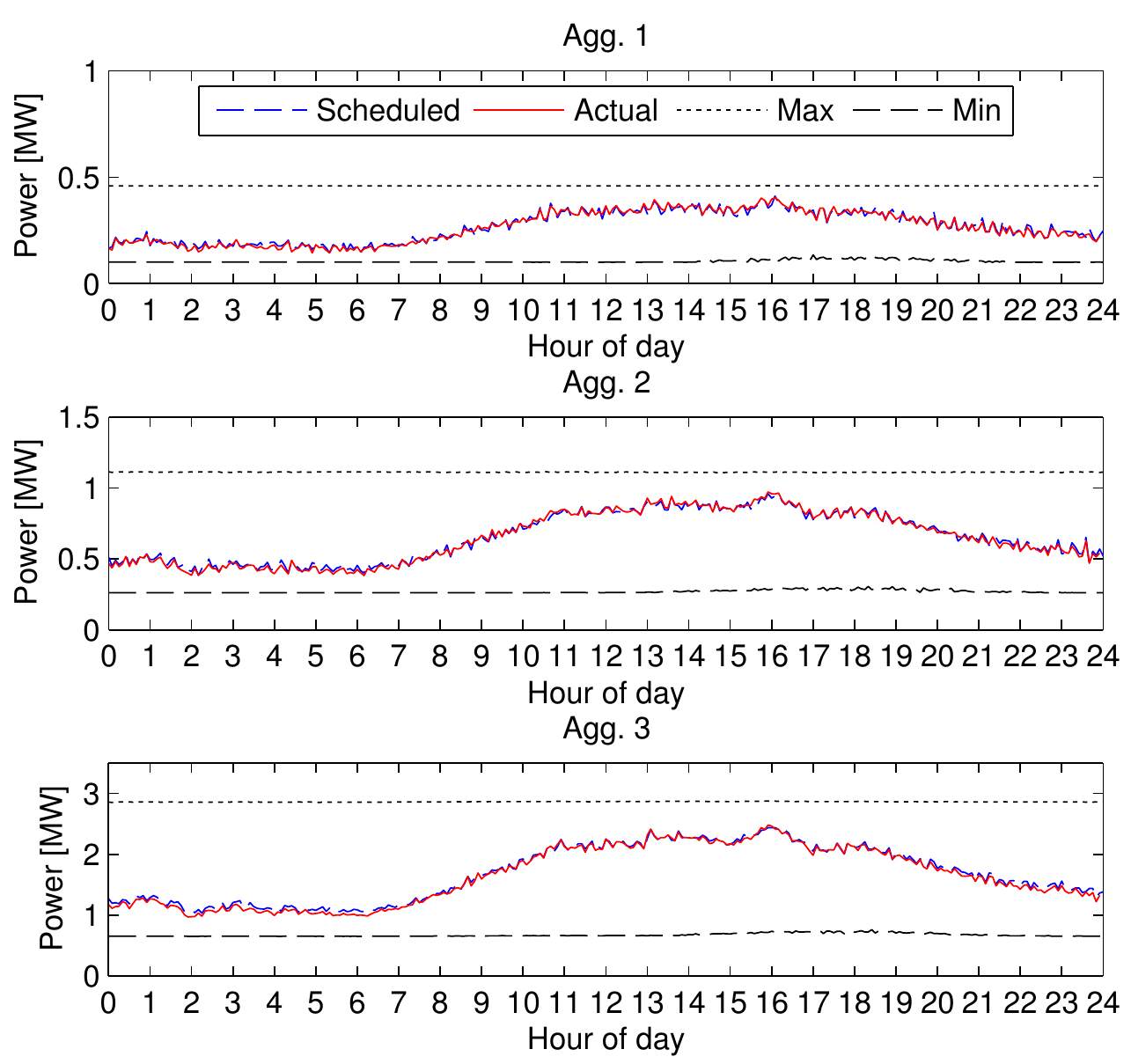}
\caption{Load under each aggregator in Case 2.}
\label{fig:Agg_load}
\end{figure}

The output of DGs is shown in Fig.~\ref{fig:DG_gen}.
DG2 is the cheapest generator and is at its maximum output almost all the time.
Other DGs generate more during peak hours, because the reference signal essentially requires more reduction from the base case during peak hours.
It can be easily verified that the marginal cost of all DGs that are not at their generation limits is the same, using the cost parameters in Table~\ref{table-DG}.
The scheduled and actual load from aggregators together with their dynamic capability (max and min) are plotted in Fig.~\ref{fig:Agg_load}.
The feasible load range for each AC in each time period depends on the current indoor air temperature, temperature setpoint and acceptable range, price information etc., and therefore varies significantly from one time period to another.
Nevertheless, the feasible load range from aggregating a large number of ACs does not vary much.
The actual average power consumption closely follows the desired value, which verifies the effectiveness of the proposed coordination and control.

\begin{figure}
\centering
\includegraphics[width=0.45\textwidth]{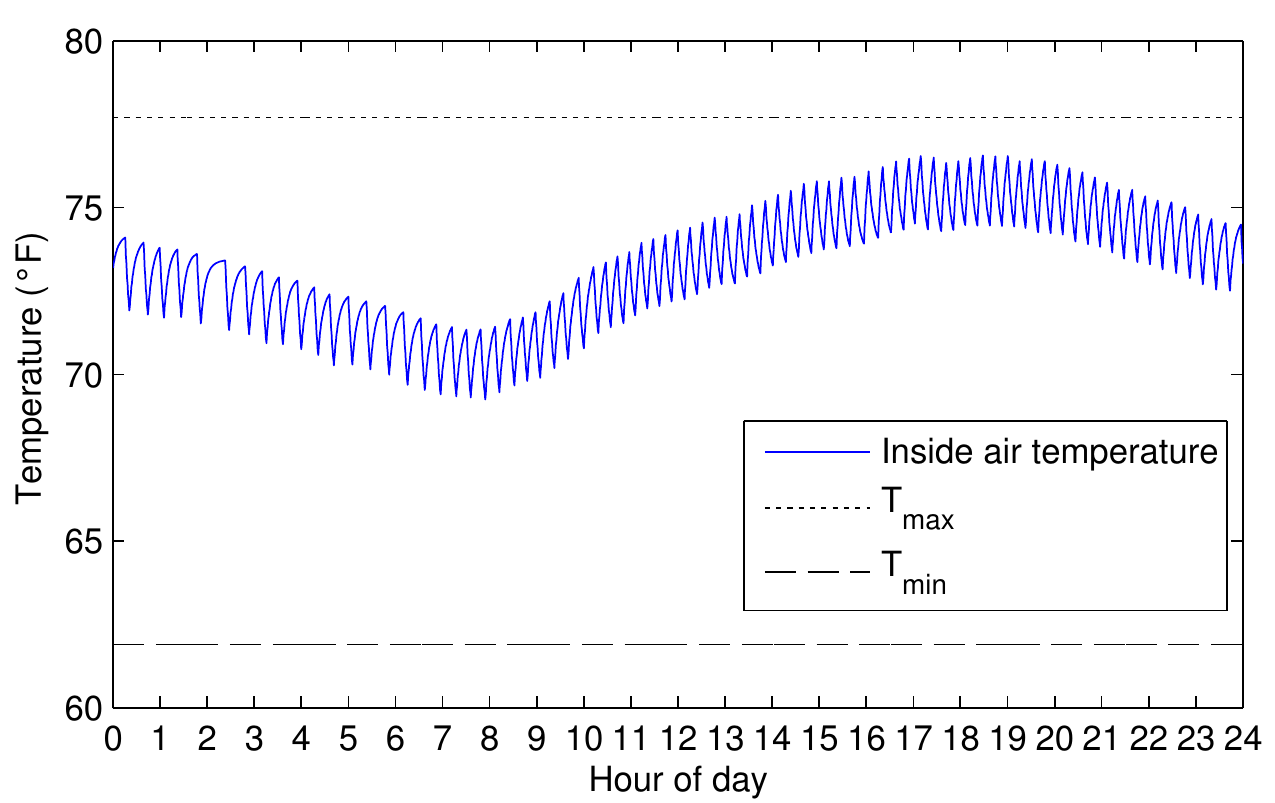}
\caption{Indoor air temperature of House~1 under Aggregator~1 in Case~2.}
\label{fig:house_Tair}
\end{figure}

The acceptable temperature settings and the simulated indoor air temperature are plotted in Fig.~\ref{fig:house_Tair} for a house under Aggregator~1.
Based on how customers value their comfort, temperature setpoint varies with the system energy cost throughout a day.
During off-peak hours when the energy price is low, the temperature setpoint and indoor air temperature are closer to the desired value, which is $72.3\,^{\circ}\mathrm{F}$ in this case.

\subsubsection{Active DER and DR Coordination with Line Capacity Constraints (Case 3)}
\begin{figure}
\centering
\includegraphics[width=0.45\textwidth]{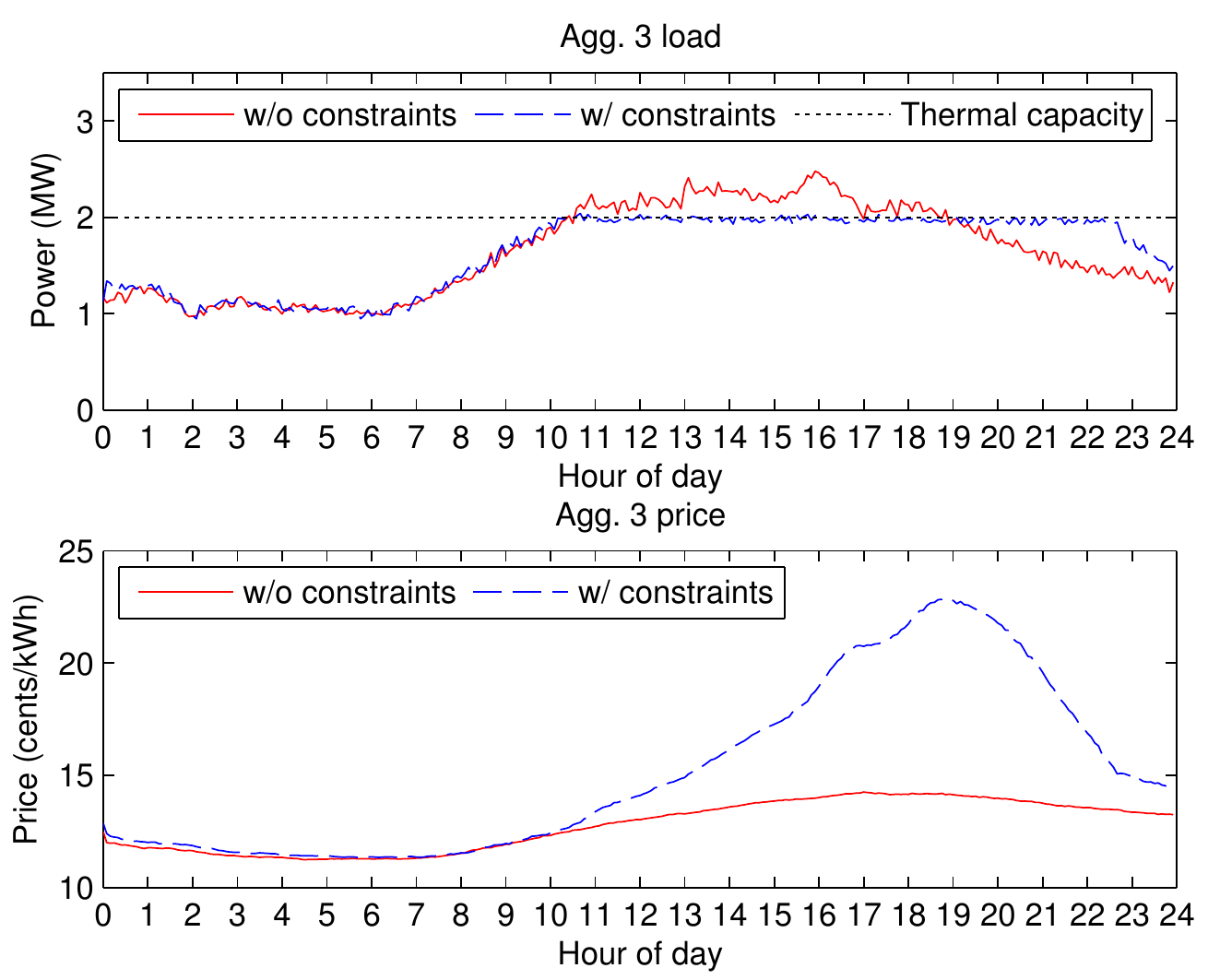}
\caption{Power consumption and clearing price under Aggregator~3 in Case~3.}
\label{fig:Agg_limits}
\end{figure}
It follows from Fig.~\ref{fig:Agg_load} that the load from Aggregator~3 exceeds 2~MW during peak hours.
Now suppose that the thermal capacity is at 2~MW for the branch that connects Aggregator~3 to the distribution system.
Since the branch connects only Aggregator~3 to the distribution system,
the power consumption by Aggregator~3 is equal to the power flow in the branch (ignoring losses for simplicity).
Therefore, the capacity limit of such a branch can be taken care of by imposing the limit to Aggregator~3's power range, i.e., modifying the maximum power consumption $p^\text{max}_i$ of Aggregator 3 in (2c).
In this case, the maximum power limit of Aggregator~3 in (2c) becomes active for some time.
The simulation results are shown in Fig.~\ref{fig:Agg_limits}.
It can be seen that the proposed framework can account for local thermal constraints as well when coordinating DERs and DR.
When congestion occurs, the active local constraint significantly raises the energy price for the load under Aggregator~3.
It should be noted that the load and price are obtained using the same optimal coordination algorithm, rather than being manually modified to meet the local constraint.
Compared with the case where there is not line capacity constraints, these load and price are just some different solutions of the optimal coordination problem (2) with some updated parameters.

\section{Conclusions}\label{sec:conclusion}
This paper presents a hierarchical control framework to integrate DR into DER coordination.
The proposed framework takes the advantage of existing coordination algorithms and device controllers, and bridges the gap between short-term scheduling and real-time control of controllable loads.
This is done by synthetically capturing the underlying economics from DR as well as detailed dynamics for device-level control.
Simulation results showed that the proposed method is capable to optimally coordinate DR with DGs and control DR in real-time to realize the desired allocation of power consumption.
The future work is to expand this framework to include distributed energy storage into this coordination between DGs and controllable loads.

\appendices
\section{Demand Curve Determination}
The demand curve of each AC as shown in Fig.~\ref{fig:ACdemand} is characterized by the maximum and minimum energy consumptions ($E_\text{max}^i$ and $E_\text{min}^i$), and the corresponding energy prices ($\lambda_\text{min}^i$ and $\lambda_\text{max}^i$), which can be calculated as follows:
\begin{itemize}
  \item $E_\text{max}^i$ and $\lambda_\text{min}^i$:
  For the $i$-th unit, the theoretical upper bound of average power consumption corresponds to the operation when the device is ON for the entire period.
      In this case, the average power consumption is simply equal to $p_c$, which is the instantaneous power when the device is ON and is constant through the 5-minute operating period.
      With $q_i(t)=1$, the closed-form analytical expression of indoor air temperature $T_a^i(t)$ can be obtained by solving~\eqref{eq:etp}.
      The setpoint $T_\text{set}$ must be low enough to satisfy \eqref{eq:on} to maintain the ON status for the entire period assuming the device is operated in cooling mode,
      \begin{equation}\label{eq:on}
        T_\text{set} \le T_\text{set}^\text{on}\,\text{,}
      \end{equation}
      where $T_\text{set}^\text{on}=\min_t\{T_a^i(t)\} - \delta/2$ when the device is initially OFF, and $T_\text{set}^\text{on}=\min_t\{T_a^i(t)\} + \delta/2$ when the device is initially ON.
        \begin{itemize}
        \item If $T_\text{set}^\text{on} \ge T_\text{min}$, where $T_\text{min}$ is the lowest acceptable temperature setpoint specified by the user, find the energy price $\lambda_\text{min}^i$ which corresponds to $T_\text{set}^\text{on}$ on the curve in Fig~\ref{fig:response}. For any price that is less than $\lambda_\text{min}^i$, AC~$i$ will be ON for the entire operating period and the maximum average power consumption $E_\text{max}^i$ is $p_c$.
        \item If $T_\text{set}^\text{on} < T_\text{min}$, the device can only be ON for part of the period, because $T_\text{set}$ ($\ge T_\text{min} >T_\text{set}^\text{on}$) will trigger the device to be OFF for at least some time during the period.
            In this case, the energy price $\lambda_\text{min}^i$ simply corresponds to $T_\text{min}$ on the curve in Fig~\ref{fig:response}. The corresponding upper bound of feasible average power consumption $E_\text{max}^i$ can be found by solving \eqref{eq:etp} and \eqref{eq:hyster} by letting $T_\text{set}=T_\text{min}$.
        \end{itemize}
    \item $E_\text{min}^i$ and $\lambda_\text{max}^i$: The theoretical lower bound of average power consumption is zero and it corresponds to the operation when the device is OFF for the entire period.
The setpoint $T_\text{set}$ must satisfy \eqref{eq:off} to maintain the OFF status for the entire period when the device is operated in cooling mode,
\begin{equation}
\label{eq:off}
  T_\text{set} \ge T_\text{set}^\text{off}\,\text{,}
\end{equation}
where $T_\text{set}^\text{off}=\max_t\{T_a^i(t)\} - \delta/2$ when the device is initially OFF, and $T_\text{set}^\text{off}=\max_t\{T_a^i(t)\} + \delta/2$ when the device is initially ON.
\begin{itemize}
  \item If $T_\text{set}^\text{off} \le T_\text{max}$, where $T_\text{max}$ is the highest acceptable temperature setpoint specified by the user,
  find the energy price $\lambda_\text{max}^i$ which corresponds to $T_\text{set}^\text{off}$ on the curve in Fig~\ref{fig:response}.
    For any price that is higher than $\lambda_\text{max}^i$, AC~$i$ will be OFF for the entire operating period and there is no energy consumption.
  \item If $T_\text{set}^\text{off} > T_\text{max}$, the device can only be OFF for part of the period, because $T_\text{set}$ ($\le T_\text{max} < T_\text{set}^\text{off}$) will trigger the device to be ON for at least some time during the period.
      In this case, the energy price $\lambda_\text{max}^i$ simply corresponds to $T_\text{max}$ on the curve in Fig~\ref{fig:response}.
      The corresponding lower bound of feasible average power consumption $E_\text{min}^i$ can be found by solving \eqref{eq:etp} and \eqref{eq:hyster} by letting $T_\text{set}=T_\text{max}$.
\end{itemize}

\item For any point between $\lambda_\text{min}^i$ and $\lambda_\text{max}^i$, the device will only be ON for part of the period.
The closed-form analytical expression does not exist, but the numerical method can be used to find the corresponding average power consumption.
This process needs to be repeated for a large number of points to accurately characterize the demand curve.
On the other hand, submitting the entire demand curve to aggregators requires sending all the points, which also burdens the communication network.
In fact, many practical applications do not require the entire demand curve, and its approximation should good enough for engineering purpose.
For example, the demand curve shown in Fig.~\ref{fig:ACdemand} can be approximated by a step curve as long as the difference between $\lambda_\text{min}^i$ and $\lambda_\text{max}^i$ is small.
In this case, only three quantities including $E_\text{min}^i$, $E_\text{max}^i$ and $(\lambda_\text{min}^i+\lambda_\text{max}^i)/2$ are required to characterize the approximated step curve,
thereby reducing communication requirement and computational cost.
\end{itemize}


\end{document}